\documentclass[11pt]{article}
\newcommand{ \newsection}[1]{ \setcounter{equation}{0} \section{ #1} }

\def\keywords{ \if@twocolumn
\section*{Keywords}
\else \small
\begin{center}
{ \bf Keywords\vspace{-.5em}\vspace{0pt}}
\end{center}
\center \fi}
\def\endkeywords{ \if@twocolumn\else\endcenter\fi}
\usepackage{epsfig}
\usepackage{graphicx}
\begin{document}
\title{\bf Non-stationary subdivision schemes originated from uniform trigonometric B-spline}
\author{Shahid S. Siddiqi  \thanks{
 Department of Mathematics, University of the Punjab, Lahore 54590,
 Pakistan. ~~~~~~~~~~~~~~~~~~~~~~~~~~~~~~~~~~~~~~~~~~~~~~~~~~~~Email: shahidsiddiqiprof@yahoo.co.uk.} and Muhammad Younis \thanks{
 Email: younis.pu@gmail.com}}
\date{}
\maketitle
\begin{center}
\begin{minipage}{5.0in}
\ \ \ \\
{\bf Keywords:} {\small
 binary, approximation, non-stationary schemes, uniform trigonometric B-spline, convergence and smoothness} \\
\begin{abstract}
\noindent

The paper proposes, an algorithm to produce novel $m$-point (for
any integer $m \geq 2$) binary non-stationary subdivision scheme.
It has been developed using uniform trigonometric B-spline basis functions
and smoothness is being analyzed using the theory of
asymptotically equivalence. The results show that the most of
well-known binary approximating schemes can be considered as the non-stationary counterpart of the
proposed algorithm.

Furthermore, the schemes developed by the proposed algorithm has the ability to reproduce
 or regenerate the conic sections, trigonometric polynomials and trigonometric splines as well.
 Some examples are considered, by choosing an appropriate tension parameter $0<\alpha<\pi/3
$, to show the usefulness.
\end{abstract}
\end{minipage}
\end{center}
\newsection{Introduction}

~~~~Subdivision scheme is one of the most important and
significant modelling tool to create smooth curves from initial
control polygon by subdividing them according to some refining
rules, recursively. These refining rules take the initial polygon
to produce a sequence of finer polygons converging to a smooth
limiting curve.
%

In the field of non-stationary subdivision schemes, Beccari
\textit{et al.}~\cite{BC077} presented a 4-point binary non-stationary
interpolating subdivision scheme, using tension parameter, that
was capable of producing certain families of conics and cubic
polynomials. They also developed a 4-point ternary
interpolating non-stationary subdivision scheme in the same year
that generate $C^2$ continuous limit curves showing considerable
variation of shapes with a tension parameter~\cite{BC07}. A new
family of 6-point interpolatory non-stationary subdivision scheme
was introduced by Conti and Romani~\cite{COS10}. It was presented
using cubic exponential B-spline symbol generating functions that
can reproduce conic sections. Conti and Romani \cite{COS11} discussed algebraic conditions on non-stationary subdivision symbols
 for exponential polynomial reproduction.

Since, non-stationary schemes have proven to be efficient
iterative algorithms to construct special classes of curves. One
of the important capability is the reproduction or regeneration of
trigonometric polynomials, trigonometric splines and conic
sections, in particular circles, ellipses etc. So, in this
article, an algorithm has been introduced to produce $m-$point
binary approximating non-stationary schemes, for any integer $m
\geq 2$, using the uniform trigonometric B-spline basis function
of order $m-1$. The proposed algorithm can be considered as the
non-stationary counterpart of the well known binary approximating
schemes introduced by Chakin \cite{GM74} and Siddiqi with his
different co-authors \cite{SSK10, SSK09, SS11, SY12, SSN08}, after
setting different values of $m$ in proposed algorithm (for details
see table 5.1). Moreover, the proposed algorithm can also be
considered as generalization form of the 2-point and 3-point
non-stationary schemes presented by Daniel and Shunmugaraj
\cite{SDPS09}.

The paper is organized as follows, in section $2$ the basic notion
and definitions of binary subdivision scheme are considered. The
algorithm, to produce $m-$point binary non-stationary scheme, is
presented in section $3$. Some example are considered, to
construct
 the masks of 2-point, 3-point and 4-point schemes, in section 4.
 The convergence and smoothness of the schemes are being calculated in
 section $5$. Some properties and advantages of proposed algorithm
are being discussed in section $6$. The conclusion is drawn in
section $7$.

\section{Preliminaries}

~~~~In univariate subdivision scheme, following the notion and definitions introduced in
\cite{ND02}, the set of control points $ \{f_{i}^{k} \in R\mid
i\in Z\}$ of polygon at $k^{th}$ level is mapped to a refined
polygon to generate the new set of control points $
\{f_{i}^{k+1} \in R\mid i\in Z\}$ at the $(k+1)^{st}$ level by
applying the following repeated application of the refinement rule
\begin{eqnarray}
f^{k+1}_{i}=\{S_{a^{(k)}}f^{k}\}_{i}=\sum_{j\in
Z}a^{(k)}_{i-2j}f^{k}_{j}&& \forall~i \in Z,
\end{eqnarray}
with the sequence of finite sets of real coefficients
$\{a_{i}^{(k)}, k\geq 0\}$ constitute the so-called mask of
subdivision scheme. If the mask of a scheme are independent of
$k$, namely if $a^{(k)}=a$ for all $k\geq0$ then it is called
stationary $ S_{a}$, otherwise it is called non-stationary
$S_{a^{(k)}}$.

The convergent subdivision scheme, formally denoted by $f^{k+1}=S_{a^{(k)}}f^k$, with the corresponding mask
$a_{i}^{(k)}~i\in Z$ necessarily satisfies
 \begin{eqnarray*}
 \sum_{j\in Z} a^{(k)}_{2j}\cong\sum_{j\in Z} a^{(k)}_{2j+1}\cong 1.
\end{eqnarray*}
A binary non-stationary subdivision scheme
 $S_{a^{(k)}}$ is said to be convergent if for every initial
 data $ f^{0}\in l^{\infty}$ there exits a continuous limit function $ g\in C^{m}(R) $
 such that
\begin{eqnarray*}
  \lim_{k\rightarrow\infty}\sup_{i\in Z}|f_{i}^{k}-g(2^{-k}i)|=0,
\end{eqnarray*}
and $ g $ is not identically zero for some initial data $
f^{0}$.

We also recall that two univariate binary schemes
 $S_{a^{(k)}}$ and $ S_{b^{(k)}}$ are said to be asymptotically
 equivalent if
 \begin{eqnarray*}
\sum^{\infty}_{k=1}\|S_{a^{(k)}}-S_{b^{(k)}}\|_{\infty}<\infty ,
\end{eqnarray*}
~~~~~~~~~~~~~~ where
\begin{eqnarray*}
\|S_{a^{(k)}}\|_{\infty}=max\{\sum_{i \in Z}|a_{2i}^{(k)}|, \sum_{i
\in Z}|a_{2i+1}^{(k)}|\}\end{eqnarray*}
 \textbf{Theorem 2.1.} \textit{The
non-stationary scheme $ S_{a^{(k)}}$ and stationary scheme $
S_{a}$  are said to be asymptotically equivalent schemes, if
they have finite masks of the same support. The stationary scheme
$S_{a}$ is $ C^{m}$ and}
 \begin{eqnarray*}
\sum^{\infty}_{k=0}2^{mk}\|S_{a^{(k)}}-S_{a}\|_{\infty}<\infty,
\end{eqnarray*}
\textit{then non-stationary scheme $ S_{a^{(k)}}$ is also said to be $
C^{m}$.}

\section{The algorithm for approximating schemes}

~~~~ In this section, an algorithm has been introduced to produce
$m$-point binary non-stationary subdivision schemes (for any
integer $m \geq 2$) which can generate the families of $C^{m-1}$
limiting curves by choosing a tension parameter $0<\alpha<\pi/3$.
The algorithm has been established using uniform trigonometric
B-splines of order n. So, in view of \textit{Koch et
al.}~\cite{KC95}, trigonometric B-splines can be defined as
follow. Let $m>n>0$ and $0<\alpha<\frac{\pi}{n}$, then Uniform
Trigonometric B-splines $\{T_{j}^{n}(x;\alpha)\}_{j=1}^{m}$ of
order $n $ associated with the knot sequence $\Delta := \{t_{i} =
i\alpha: i=0, 1, 2, ..., m+n\}$ with the mesh size $\alpha$ are
defined by the recurrence relation,
 \textbf{}\\
\begin{eqnarray*}
T^{1}_{0}(x;\alpha)= \left\{\begin{array}{rr}
1, &  x\in[0,\alpha)\\
0, &  otherwise,\\
\end{array}\right.
\end{eqnarray*}
for $1<r\leq n,$
\begin{eqnarray}
T^{r}_{0}(x;\alpha)=
\frac{1}{\sin((r-1)\alpha)}\{\sin(x)T_{0}^{r-1}(x;\alpha)+\sin(t_{r}-x)T_{0}^{r-1}(x-\alpha;\alpha)\}
\end{eqnarray}
and $T^{r}_{j}(x;\alpha)=T^{r}_{0}(x-j\alpha;\alpha),$ for $j=1,
2, ...,m$. The trigonometric B-spline $T^{n}_{j}(x;\alpha)$ is
supported on $[t_{j},t_{j+n}]$ and it is the interior of its
support. Moreover, $\{T^{n}_{j}\}_{j=1}^{m}$  are linearly
independent set on the interval $[t_{n-1},t_{m+1}]$. Hence, on
this interval, any uniform trigonometric spline $S(x)$ has a
unique representation of the form
$S(x)=\sum_{j=0}^{m}p_{j}T_{j}^{n}(x;\alpha), ~~p_{j}\in R$.

To obtain the mask $a_{i}^{k}(\alpha)=a_{i}^{k},~~
i=0,1,2,...,m-1$ we use the following recreance relation, for any
value of $k$,
\begin{eqnarray}
a_{i}^{k}(\alpha)=T_{0}^{m}\left(
(m-i-1)\frac{\alpha}{2^{k}}+\frac{\alpha}{2^{k+2}};\frac{\alpha}{2^{k}}\right)~~~~~i=0, 1, 2, ...,m-1
\end{eqnarray}
where $T_{0}^{m}\left(x;\frac{\alpha}{2^{k}}\right)$, with mesh
size $(\frac{\alpha}{2^{k}})$, is a trigonometric B-spline basis
function of order $m-1$ and can be calculated form equation (3.2).
It can also be observed that the schemes produced by the proposed
algorithm do not have the convex hull and affine invariance
properties. Since the sums of the weights of the obtained schemes
at $kth$ level are not equal to unity. To get sum of the mask equal to unity, the corresponding normalized scheme can be obtained (see \cite{SDPS09}). In the following, some
examples are considered to produce the masks of 2-point, 3-point
and 4-point binary approximating schemes after setting $m=2$, 3
and 4, respectively, in above recurrence relation.

\section{Construction of the schemes}

In this section, some applications are considered to construct the masks of
2-point, 3-point and 4-point approximating non-stationary schemes.

\subsection{The 2-point Approximating Scheme}

The linear trigonometric B-spline basis function $T_{0}^{2}\left(
x;\frac{\alpha}{2^{k}}\right)$, with mesh size
$(\frac{\alpha}{2^{k}})$, can be calculated by setting $m=2$ in
relation (3.3). The 2-point binary non-stationary scheme (which is
also called corner cutting scheme) with mask
$a_{i}^{(k)}(\alpha)=a_{i}^{(k)}, i=0,1$  is defined,
for any value of $k$ and $\alpha \in ]0, \frac{\pi}{3}[$, as

\begin{eqnarray}
\left.\begin{array}{rrr}
f^{k+1}_{2i}&=&a_{0}^{(k)}f^{k}_{i}+a_{1}^{(k)}f^{k}_{i+1}\\
f^{k+1}_{2i+1}&=&a_{1}^{(k)}f^{k}_{i}+a_{0}^{(k)}f^{k}_{i+1}
\end{array}\right\}\
\end{eqnarray}
where
\begin{eqnarray*}
a^{(k)}_{0}=
\frac{\sin\frac{{3\alpha}}{2^{k+2}}}{\sin\frac{\alpha}{2^{k}}},~~~~
a^{(k)}_{1}=
\frac{\sin\frac{{\alpha}}{2^{k+2}}}{\sin\frac{\alpha}{2^{k}}}
\end{eqnarray*}
\textbf{Theorem 4.1.1.} \textit{ The 2-point binary non-stationary scheme defined above converges and has smoothness, $C^1$ for the range $\alpha\in ]0,\frac{\pi}{3}[$. }\\
\textbf{Proof.} see~\cite{SDPS09}.\\\\
\textbf{Remark 4.1.2.} It can be observed that the mask of above 2-point normalized scheme converges to the mask of the famous corner cutting scheme introduced by Chaikin~\cite{GM74}. Moreover,
it can also be observed that the mask of binary 2-point non-stationary scheme of Daniel and Shummugaraj~\cite{SDPS09} can be calculated, after setting $m=2$ in (4). Hence the proposed
scheme can be considered as the generalized form of non-stationary scheme presented by Daniel and Shummugaraj and non-stationary counterpart of famous Chaikin's scheme~\cite{GM74}.

\subsection{The 3-point Approximating Scheme}

To get quadratic trigonometric B-spline basis function,
$T_{0}^{3}\left(x;\frac{\alpha}{2^{k}}\right)$ with mesh size
$(\frac{\alpha}{2^{k}})$, we take $m=3$ in equation (3.3). The masks
$a_{i}^{(k)}(\alpha)=a_{i}^{(k)}, i=0,1,2$ of the
proposed binary 3-point scheme can be calculated from quadratic
trigonometric B-spline function.

The 3-point non-stationary scheme is defined, for some value of
$\alpha \in ]0, \frac{\pi}{3}[$, as follow
\begin{eqnarray}
\left.\begin{array}{rrr}
f^{k+1}_{2i}&=&a_{0}^{(k)}f^{k}_{i-1}+a_{1}^{(k)}f^{k}_{i}+a_{2}^{(k)}f^{k}_{i+1}\\
f^{k+1}_{2i+1}&=&a_{2}^{(k)}f^{k}_{i-i}+a_{1}^{(k)}f^{k}_{i}+a_{0}^{(k)}f^{k}_{i+1}
\end{array}\right\}\
\end{eqnarray}
where
\begin{eqnarray*}
a^{(k)}_{0}=
\frac{\sin^{2}\frac{{3\alpha}}{2^{k+2}}}{\sin\frac{\alpha}{2^{k}}\sin\frac{2\alpha}{2^{k}}}~~~~
a^{(k)}_{1}=\frac{\sin\frac{{3\alpha}}{2^{k+2}}\sin\frac{{5\alpha}}{2^{k+2}}+
\sin\frac{{\alpha}}{2^{k+2}}\sin\frac{{7\alpha}}{2^{k+2}}}{\sin\frac{\alpha}{2^{k}}\sin\frac{2\alpha}{2^{k}}}~~~~
a^{(k)}_{2}=
\frac{\sin^{2}\frac{{\alpha}}{2^{k+2}}}{\sin\frac{\alpha}{2^{k}}\sin\frac{2\alpha}{2^{k}}}~~~~
\end{eqnarray*}
\textbf{Theorem 4.1.3.} \textit{ The 3-point binary non-stationary scheme defined above converges and has smoothness, $C^2$ for the range $\alpha\in ]0,\frac{\pi}{3}[$. }\\
\textbf{Proof.} see~\cite{SDPS09}.\\\\
\textbf{Remark 4.1.4.} The proposed scheme (4.5) is considered as the generalized
form of the non-stationary 3-point scheme developed by Daniel and Shummugaraj in~\cite{SDPS09}.
Furthermore, the proposed 3-point scheme (4.5) can be considered as
the non-stationary counterpart of the stationary scheme \cite{SS11}.

\subsection{The 4-point approximating scheme}

In this section, a 4-point binary approximating non-stationary subdivision scheme
is presented and masks of the proposed 4-point binary
scheme can be calculated, for any value of $k$, using the relation
(4). Where $T_{0}^{4}\left( x;\frac{\alpha}{2^{k}}\right)$, with
mesh size $(\frac{\alpha}{2^{k}})$, is called the cubic trigonometric
B-spline basis functions and can be calculated from the recurrence
relation (3). The proposed scheme is defined, for some value of
$\alpha \in ]0, \frac{\pi}{3}[$, as
\begin{eqnarray}
\left.\begin{array}{rrr}
f^{k+1}_{2i}&=&a_{0}^{(k)}f^{k}_{i-1}+a_{1}^{(k)}f^{k}_{i}+a_{2}^{(k)}f^{k}_{i+1}+a_{3}^{(k)}f^{k}_{i+2}\\
f^{k+1}_{2i+1}&=&a_{3}^{(k)}f^{k}_{i-i}+a_{2}^{(k)}f^{k}_{i}+a_{1}^{(k)}f^{k}_{i+1}+a_{0}^{(k)}f^{k}_{i+2}
\end{array}\right\}\
\end{eqnarray}
where
\begin{eqnarray*}
a^{(k)}_{0}&=&
\frac{\sin^{3}\frac{{\alpha}}{2^{k+2}}}{\sin\frac{\alpha}{2^{k}}\sin\frac{2\alpha}{2^{k}}\sin\frac{3\alpha}{2^{k}}}\\
a^{(k)}_{1}&=&\frac{\sin\frac{{3\alpha}}{2^{k+2}}\sin^{2}\frac{{5\alpha}}{2^{k+2}}+
\sin\frac{{\alpha}}{2^{k+2}}\sin\frac{{5\alpha}}{2^{k+2}}\sin\frac{{7\alpha}}{2^{k+2}}+
\sin^{2}\frac{{\alpha}}{2^{k+2}}\sin\frac{{11\alpha}}{2^{k+2}}}{\sin\frac{\alpha}{2^{k}}\sin\frac{2\alpha}{2^{k}}\sin\frac{3\alpha}{2^{k}}}\\
a^{(k)}_{2}&=&
\frac{\sin^{2}\frac{{3\alpha}}{2^{k+2}}\sin\frac{{9\alpha}}{2^{k+2}}+
\sin\frac{{3\alpha}}{2^{k+2}}\sin\frac{{5\alpha}}{2^{k+2}}\sin\frac{{7\alpha}}{2^{k+2}}+
\sin\frac{{\alpha}}{2^{k+2}}\sin^{2}\frac{{7\alpha}}{2^{k+2}}}{\sin\frac{\alpha}{2^{k}}\sin\frac{2\alpha}{2^{k}}\sin\frac{3\alpha}{2^{k}}}\\
a^{(k)}_{3}&=&
\frac{\sin^{3}\frac{{3\alpha}}{2^{k+2}}}{\sin\frac{\alpha}{2^{k}}\sin\frac{2\alpha}{2^{k}}\sin\frac{3\alpha}{2^{k}}}
\end{eqnarray*}
the proposed 4-point scheme can also be considered as the general
form of stationary 4-point binary approximating scheme, which was
introduced by Siddiqi and Ahmad \cite{SSK10}. The subdivision rules
to refine the control polygon are defined as
\begin{eqnarray}
\left.\begin{array}{rrr}
f^{k+1}_{2i}&=&\frac{1}{384}f^{k}_{i-1}+\frac{121}{384}f^{k}_{i}+\frac{235}{384}f^{k}_{i+1}+\frac{27}{384}f^{k}_{i+2}\\
f^{k+1}_{2i+1}&=&\frac{27}{384}p^{k}_{i-i}+\frac{235}{384}f^{k}_{i}+\frac{121}{384}f^{k}_{i+1}+\frac{1}{384}f^{k}_{i+2}
\end{array}\right\}\
\end{eqnarray}
As the weights of the mask of the proposed scheme (4.6) are bounded
by the coefficient of the mask of the above scheme (4.7). So, we can write as
\begin{eqnarray*}
a_{0}^{(k)}\rightarrow \frac{1}{384},~~
a_{1}^{(k)}\rightarrow \frac{121}{384},
~~a_{2}^{(k)}\rightarrow
\frac{235}{384},~~a_{3}^{(k)}\rightarrow \frac{27}{384},~~
\end{eqnarray*}

The proofs of $a_{0}^{(k)}\rightarrow \frac{1}{384}$, $a_{1}^{(k)}\rightarrow
\frac{121 }{384}$, $a_{2}^{(k)}\rightarrow \frac{235}{384}$ and
$a_{3}^{(k)}\rightarrow \frac{27}{384}$ can be followed
from the lemma (5.1.1).

\section{Convergence Analysis}

The theory of asymptotic
equivalence is used to investigate the convergence and smoothness
of the proposed scheme following~\cite{ND95}. Some estimations of $\gamma_{i}^{k}, i=0,1,2,
...,m-1$ are bring into play to prove the convergence of the proposed schemes. To establish the estimations some inequalities are being considered.

\begin{eqnarray*}
\frac{\sin a}{\sin b}\geq\frac{a}{b} &for&0<a<b<\frac{\pi}{2}\\
\theta \csc\theta \leq t\csc t &for& 0<\theta<t<\frac{\pi}{2}
\end{eqnarray*}
~~~~~~~~~~~~~~~~~~~~~~~and
\begin{eqnarray*}
\cos x\leq\frac{\sin x}{x}   &for& 0<x<\frac{\pi}{2}
\end{eqnarray*}

\subsection{Convergence Analysis of 4-point Scheme}
To prove the convergence and smoothness of scheme (4.6),
estimations of $a_{i}^{(k)}, i=0,1,2,3 $
are being calculated in the following lemmas.\\\\
\textbf{Lemma 5.1.1.} \textit{ For $k\geq 0$} and
\textit{$0<\alpha<\frac{\pi}{2}$}
\begin{eqnarray*} \textbf{(i)}~~
\frac{1}{384}\leq a_{0}^{(k)}\leq\frac{1}{384}\frac{1}{\cos^{3}(\frac{6\alpha}{2^{k}})}~~~~~~~~
\textbf{(ii)}~~
\frac{121}{384}\leq a_{1}^{(k)}\leq\frac{121}{384}\frac{1}{\cos^{3}(\frac{6\alpha}{2^{k}})}\\
\textbf{(iii)}~~
\frac{235}{384}\leq a_{2}^{(k)}\leq\frac{235}{384}\frac{1}{\cos^{3}(\frac{6\alpha}{2^{k}})}~~~~~~
\textbf{(iv)}~~
\frac{27}{384}\leq a_{3}^{(k)}\leq\frac{27}{384}\frac{1}{\cos^{3}(\frac{6\alpha}{2^{k}})}
\end{eqnarray*}
\textbf{Proof}. To prove the inequality \textbf{(i)}
\begin{eqnarray*}
a^{(k)}_{0}&=&\frac{\sin^{3}(\frac{\alpha}{2^{k+3}})}
{\sin(\frac{\alpha}{2^{k}})\sin(\frac{2\alpha}{2^{k}})\sin(\frac{3\alpha}{2^{k}})}
\geq \frac{(\frac{\alpha}{2^{k+3}})^{3}}
{(\frac{\alpha}{2^{k}})(\frac{2\alpha}{2^{k}})(\frac{3\alpha}{2^{k}})}
=\frac{1}{384}\\
\end{eqnarray*}
and
\begin{eqnarray*}
a^{(k)}_{0}&\leq&\frac{\alpha^{3}}{2^{3k+6}}\csc\left(\frac{\alpha}{2^{k}}\right)
\csc\left(\frac{2\alpha}{2^{k}}\right)\csc\left(\frac{3\alpha}{2^{k}}\right)
\leq\frac{\alpha^{3}}{2^{3k+6}}~~36\csc^{3}\left(\frac{6\alpha}{2^{k}}\right)\\
&\leq&\frac{\alpha^{3}}{2^{3k+6}}~~36~~\frac{1}{(\frac{6\alpha}{2^{k}})^{3}\cos^{3}(\frac{6\alpha}{2^{k}})}
\leq\frac{1}{384}\frac{1}{\cos^{3}(\frac{6\alpha}{2^{k}})}
\end{eqnarray*}
The proofs of \textbf{(ii), (iii)} and \textbf{(iv)} can be
obtained similarly.\\\\
\textbf{Lemma 5.1.2.} For some constants $ C_{0},C_{1},C_{2}$ and
$C_{3}$ independent of $k$, we have
\begin{eqnarray*}
\textbf{(i)} ~~ |a^{(k)}_{0}-\frac{1}{384}|\leq
C_{0}\frac{1}{2^{3k}}~~~~~~~~
\textbf{(ii)} ~~ |a^{(k)}_{1}-\frac{121}{384}|\leq C_{1}\frac{1}{2^{3k}}\\
\textbf{(iii)} ~~ |a^{(k)}_{2}-\frac{235}{384}|\leq
C_{2}\frac{1}{2^{3k}}~~~~~~
 \textbf{(iv)} ~~
|a^{(k)}_{3}-\frac{27}{384}|\leq C_{3}\frac{1}{2^{3k}}
\end{eqnarray*}
\textbf{Proof}. To prove the inequality \textbf{(i)} use Lemma
(5.1.1),
 \begin{eqnarray*}
\textbf{(i)} ~~ |a^{(k)}_{0}-\frac{1}{384}|&\leq&
\frac{1}{384}\left(\frac{1-\cos^{3}(\frac{6\alpha}{2^k})}{\cos^{3}\alpha}\right) \\
&\leq&\frac{1}{384}\left(\frac{3\sin^{3}(\frac{6\alpha}{2^k})}{2\cos^{3}{\alpha}}\right)\\
&\leq&\frac{9\alpha^{2}}{64\cos^{3}(\alpha)}\frac{1}{2^{2k}}
\end{eqnarray*}
The proofs of \textbf{(ii), (iii)} and \textbf{(iv)} can be obtained similarly.\\\\
\textbf{Lemma 5.1.3.} \textit{ The laurent polynomial $b_{1}^{(k)}(z)$
of the scheme $\{S_{b^{(k)}}\}$ at the $kth$ level can be written as
$ b_{1}^{(k)}(z)=\left(\frac{1+z}{2}\right)b^{(k)}(z)$}, where
\begin{eqnarray*}
b^{(k)}(z)=2&\{&a^{(k)}_{0}z^{-2}+(a^{(k)}_{3}-a^{(k)}_{0})z^{-1}+(a^{(k)}_{0}+a^{(k)}_{1}-a^{k}_{3})
+(-a^{(k)}_{0}-a^{(k)}_{1}+a^{(k)}_{2}+a^{k}_{3})z\\
&+&(a^{(k)}_{0}+a^{(k)}_{1}-a^{(k)}_{3})z^{2}+(a^{(k)}_{3}-a^{(k)}_{0})z^{3}+
a^{(k)}_{0}z^{4}\}
\end{eqnarray*}
\textbf{Proof}. Since,
\begin{eqnarray*}
b_{1}^{(k)}(z)=a^{(k)}_{3}z^{-4}+a^{(k)}_{0}z^{-3}+a^{(k)}_{2}z^{-2}+a^{(k)}_{1}z^{-1}
+a^{(k)}_{1}+a^{(k)}_{2}z+a^{(k)}_{0}z^{2}
+a^{(k)}_{3}z^{3}
\end{eqnarray*}
Therefore using
$a_{0}^{(k)}+a_{1}^{(k)}+a_{2}^{(k)}+a_{3}^{(k)}=1$, it can be proved.\\\\
\textbf{Lemma 5.1.4.} \textit{ The laurent polynomial $b_{1}(z)$ of
the scheme $\{S_{b}\}$ at the $kth$ level can be written as $
b_{1}(z)=\left(\frac{1+z}{2}\right)b(z)$}, where
\begin{eqnarray*}
b(z)=\frac{1}{192}\{z^{-2}+26z^{-1}+95 +140z
+95z^{2}+26z^{3}+z^{4}\}
\end{eqnarray*}
\textbf{Proof}. To prove that the subdivision scheme $\{S_{b}\}$
corresponding to the symbol $b(z)$ is $C^{3}$, we have
\begin{eqnarray*}
d(z)&=&\frac{2b(z)}{(1+z)^{4}}\\
 &=&\frac{1}{48}\{z^{-2}+22z^{-1}+1\}
\end{eqnarray*}
Since the norm of the subdivision scheme $\{S_{d}\}$ is
\begin{eqnarray*}
\|S_{d}\|_{\infty}&=&max\{\sum_{i \in Z}|d_{2i}^{k}|, \sum_{i \in
Z}|d_{2i+1}^{k}|\}\\
&=&max\{\frac{1}{12},\frac{11}{12}\}<1
\end{eqnarray*}
So in view of Dyn \cite{ND02}, the stationary scheme is $C^{3}$.\\\\
\textbf{Theorem 5.1.5.} \textit{ The 4-point non-stationary scheme defined in Eq. (4.6) converges and has smoothness, $C^4$ for the range $\alpha\in ]0,\frac{\pi}{3}[$. }\\\\
\textbf {Proof.} To prove the proposed scheme to be $ C^{4}$, it
is sufficient to show that the scheme corresponding to the symbol $b^{k}(z)$
is $ C^{3}$ (in view of Theorem $(8)$ given by Dyn and Levin \cite{ND95}.
  Since $ \{S_{b}\}$ is $ C^{3}$ by
Lemma (5.1.4). So, it is sufficient to show for the convergence of
binary non-stationary scheme (4.6) that,
 \begin{eqnarray*}
\sum^{\infty}_{k=0}2^{3k}\|S_{b^{(k)}}-S_{b}\|_{\infty}<\infty
\end{eqnarray*}
~~~~~~~~~~~Where,
 \begin{eqnarray*}
\|S_{b^{(k)}}-S_{b}\|_{\infty}=max\{\sum^{\infty}_{j=Z}|b^{k}_{i+2j}-b_{i+2j}|:i=0,1\}
\end{eqnarray*}
Following Lemmas $(5.1.3)$ and $(5.1.4)$, we have
 \begin{eqnarray*}
\sum^{\infty}_{j\in
Z}\left|b^{(k)}_{2j}-b_{2j}\right|&=&4\left|a^{(k)}_{0}-\frac{1}{384}\right|+4\left|a^{(k)}_{0}+a^{(k)}_{1}-a^{(k)}_{3}-\frac{95}{384}\right|\\
&=&8\left|a^{(k)}_{0}-\frac{1}{384}\right|+4\left|a^{(k)}_{1}-\frac{121}{384}\right|+4\left|a^{(k)}_{3}-\frac{27}{384}\right|
\end{eqnarray*}
and similarly, it may be noted that
 \begin{eqnarray*}
\sum^{\infty}_{j\in Z}\left|b^{k}_{2j+1}-b_{2j+1}\right|=
6\left|a^{(k)}_{0}-\frac{1}{384}\right|+2\left|a^{(k)}_{1}-\frac{121}{384}\right|+2\left|a^{(k)}_{2}-\frac{235}{384}\right|+6\left|a^{(k)}_{3}-\frac{27}{384}\right|
\end{eqnarray*}
From (i), (ii), (iii) and (iv) of lemma (5.1.2),
$\sum^{\infty}_{k=0}2^{3k}|a^{(k)}_{0}-\frac{1}{384}|<\infty$,
$\sum^{\infty}_{k=0}2^{3k}|a^{(k)}_{1}-\frac{121}{384}|<\infty$,
$\sum^{\infty}_{k=0}2^{3k}|a^{(k)}_{2}-\frac{235}{384}|<\infty$
and
$\sum^{\infty}_{k=0}2^{3k}|a^{(k)}_{3}-\frac{27}{384}|<\infty$.
Hence, it can be written as
\begin{eqnarray*}
\sum^{\infty}_{k=0}2^{3k}\|S_{b^{(k)}}-S_{b}\|_{\infty}<\infty
\end{eqnarray*}

~~~~Thus by the Theorem 2.1, $c^{(k)}(z)$ is $C^{3}$ as the
associated scheme $c(z)$ is $C^{3}$. Hence, proposed scheme is $C^{4}$.\\

\section{Properties and advantages of algorithm}

In this section, some properties like unit circle reproduction
property, symmetry of basis function and some other advantages of
the proposed algorithm are being considered.
\subsection{Reproduction of unit circle}
~~~~It can be observed that certain functions like $\cos(\alpha.)$ and $\sin(\alpha .)$ can be reproduced
by the proposed non-stationary schemes. In particular, if a set of equidistant
point $f_{i}^{0}=\cos((i-\frac{1}{2})\alpha)$ and $\alpha=\frac{2\pi}{n}$, then the limit curve
is the unit circle (see also figure 2). Similarly, $\sin((i-\frac{1}{2})\alpha)$ reproduces $\sin(\alpha.)$.\\\\
\textbf{Proposition 6.1.} The limit curves of the scheme (4.4)
reproduces the functions $\cos(\alpha.)$ and $\sin(\alpha .)$
for the data points $f_{i}^{k}=\cos((i-\frac{1}{2})\frac{\alpha}{2^k})$ and $f_{i}^{k}=\sin((i-\frac{1}{2})\frac{\alpha}{2^k})$, respectively.\\
In other wards we have to show that for $k\geq0$\\\\
(i)  $f_{i}^{k}=\cos((i-\frac{1}{2})\frac{\alpha}{2^k})$, we have for $0\leq i \leq 2^kn$ as
 \begin{eqnarray*}
f_{2i}^{k+1}=\cos\left(2i\frac{\alpha}{2^{k+1}}\right)~~~~and ~~~~f_{2i}^{k+1}=\cos\left((2i+1)\frac{\alpha}{2^{k+1}}\right)
\end{eqnarray*}
(ii) Similarly, for $f_{i}^{k}=\sin((i-\frac{1}{2})\frac{\alpha}{2^k})$, we have for $0\leq i \leq 2^kn$ as
 \begin{eqnarray*}
f_{2i}^{k+1}=\sin\left(2i\frac{\alpha}{2^{k+1}}\right)~~~~and ~~~~f_{2i}^{k+1}=\sin\left((2i+1)\frac{\alpha}{2^{k+1}}\right)
\end{eqnarray*}
\textbf{Proof:}
For any initial data of the form $f_{i}^{0}=\cos((i-\frac{1}{2})\alpha)$, it can be followed
\begin{eqnarray*}
 f_{2i}^{1}&=&a_{0}^{(0)}f_{i}^{0}+a_{1}^{(0)}f_{i+1}^{0}\\
&=&\frac{\sin(\frac{3\alpha}{4})}{\sin\alpha}\cos((i-\frac{1}{2})\alpha)+\frac{\sin(\frac{\alpha}{4})}{\sin\alpha}\cos((i+\frac{1}{2})\alpha)=\cos\left(i\alpha\right)
\end{eqnarray*}
Similarly, it can be followed for $k^{th}$ refinement
\begin{eqnarray*}
 f_{2i}^{k}=\cos\left((2i)\frac{\alpha}{2^{k+1}}\right)
\end{eqnarray*}
Analogously, it can also be proved.
\begin{eqnarray*}
 f_{2i+1}^{k}=\cos\left((2i+1)\frac{\alpha}{2^{k+1}}\right)
\end{eqnarray*}
Proof of part \textbf{(ii)} is similar. So, the scheme (4.4) can
reproduce $\sin(\alpha.)$ for the initial data
$f_{i}^{0}=\sin((i-\frac{1}{2})\alpha)$.\\\\
\textbf{Proposition 6.2.} The limit curves of the scheme (4.5)
reproduces the functions $\cos(\alpha.)$ and $\sin(\alpha .)$
for the data points $f_{i}^{k}=\cos((i+\frac{1}{2})\frac{\alpha}{2^k})$ and $f_{i}^{k}=\sin((i+\frac{1}{2})\frac{\alpha}{2^k})$, respectively.\\
It can be proceed on same way. If a set of equidistant point
$f_{i}^{0}=\cos((i-\frac{1}{2})\alpha)$ and
$\alpha=\frac{2\pi}{n}$ can be chosen by the scheme (4.5), then
the limit curve is the unit circle. Similarly, $\sin((i-\frac{1}{2})\alpha)$ reproduces $\sin(\alpha.)$.\\\\
\subsection{Symmetry of Basis Limit Function}
~~~~The basis limit
function of the scheme is the limit function for the data
 \textbf{}\\
\begin{eqnarray*}
f^{0}_{i}= \left\{\begin{array}{rr}
1, &  i=0,\\
0, &  i\neq0\\
\end{array}\right.
\end{eqnarray*}

~~In order to prove that the basis limit function is symmetric about the
Y-axis.\\\\
\textbf{Theorem 6.1} \textit{The basis limit function F
is symmetric about the Y-axis}.\\
\textbf{Proof.} The symmetry of basis limit function can be followed on the same pattern following \cite{SY12}.

\subsection{Special cases} It can be observed that the binary
subdivision schemes presented in~\cite{GM74, SDPS09, SSK10,
SSK09, SS11, SY12, SSN08} are either the special cases or can be considered the non-stationary counterpart of the stationary schemes. (See also Table 1 and Table 2).

\begin{itemize}

 \item The $C^1$ limiting curves can be obtained after taking $m=2$ in proposed algorithm.
 The obtained curves of scheme (5), taking $\alpha=\frac{\pi}{180}$, coincide with
 the limit curves of the famous corner cutting scheme of Chaikin~\cite{GM74}.

  \item After setting $m=2$ and $m=3$ in proposed algorithm, the mask of 2-point and 3-point approximating non-stationary schemes, developed by Daniel and Shunmugaraj~\cite{SDPS09}, can be obtained.

  \item The limit curves of 3-point stationary approximating
  scheme, introduced by Siddiqi and Ahamd~\cite{SS11}, coincide with the limit curves of proposed 3-point scheme, after setting $m=3$ and $\alpha=\frac{\pi}{180}$.

  \item The limit curves of proposed 4-point scheme coincide with the limit
  curves obtained by the scheme ~\cite{SSK10}, for $m=4$ and $\alpha=\frac{\pi}{180}$.

  \item The limit curves of 5-point stationary approximating
  scheme introduced in~\cite{SSN08} matched with the limit curves of proposed scheme, for setting $m=5$
  and $\alpha=\frac{\pi}{180}$.

\item The limit curves of 6-point stationary approximating
  scheme introduced in~\cite{SSK09} matched with the limit curves of proposed scheme, for setting $m=6$
   and $\alpha=\frac{\pi}{180}$.

\item The proposed algorithm can also be considered as the non-stationary counterpart of $m$-point scheme developed by Siddiqi and Younis \cite{SY12}.
\end{itemize}

\begin{table}[bht]
\caption{The proposed algorithm can be considered as the non-stationary counterpart of the following existing binary approximating schemes, for different values of $m$.} \centering
\begin{small}
\begin{tabular}{ccccc}
\hline Setting $m$ &Scheme&Type& Continuity & Counterpart\\

\hline
2&2-point&Stationary&$C^{1}$&scheme \cite{GM74}  \\
3&3-point&Stationary&$C^{2}$&scheme \cite{SS11}  \\
4&4-point&Stationary&$C^{4}$&scheme \cite{SSK10} \\
5&5-point&Stationary&$C^{4}$&scheme \cite{SSN08}\\
6&6-point&Stationary&$C^{6}$&scheme \cite{SSK09} \\
$m$&$m$-point&Stationary&$C^{m-1}$&scheme \cite{SY12} \\
\hline
\end{tabular}\end{small}
\end{table}

\begin{table}[bht]
\caption{The following existing binary non-stationary approximating schemes can be considered the special case of proposed algorithm, for different values of $m$.} \centering
\begin{small}
\begin{tabular}{ccccc}
\hline Setting $m$ &Scheme&Type &Continuity & Coincides with\\

\hline
2&2-point&Non-stationary&$C^{1}$&scheme \cite{SDPS09} \\
3&3-point&Non-stationary&$C^{2}$&scheme \cite{SDPS09} \\
\hline
\end{tabular}\end{small}
\end{table}

\section{Conclusion}

~~~~ An algorithm of $m-$point binary approximating non-stationary
subdivision scheme (for any integer $m\geq 2$) has been developed
which generates the family of $C^{m-1}$ limiting curve, for
$0<\alpha<\pi/3 $. The construction of the algorithm is associated
with trigonometric B-spline basis function. It is also evident
from the examples that the limit curves of the proposed schemes
coincide with the schemes presented in~\cite{GM74, SDPS09, SS11,
SSK10, SSK09, SSN08}. So, the proposed algorithm, for different
values of $m$, can be considered as the generalized form of the
scheme \cite{SDPS09} and non-stationary counterpart of the
stationary schemes~\cite{GM74, SSK10, SSK09, SS11, SSN08}.
Moreover, the schemes produced by the algorithm can reproduce or
regenerate the trigonometric polynomials, trigonometric splines
and conic sections as well.

\begin{thebibliography}{40}


\bibitem{BC077} C. Beccari, G. Casciola and L. Romani, A
non-stationary uniform tension controlled interpolating 4-point
scheme reproducing conics, \emph{Comput. Aided. Geom. D.} , \textbf{24(1)}
(2007), 1--9.
%
\bibitem{BC07} C. Beccari, G. Casciola and L. Romani, An
interpolating 4-point $C^2$ ternary non- stationary subdivision
scheme with tension control, \emph{Comput. Aided. Geom. D.} , \textbf{24(4)}
(2007), 210--219.

\bibitem{GM74} G.M. Chaikin, An algorithm for high speed curve generation,
 \emph{Comput. Vision. Graph.}, \textbf{3(4)} (1974), 346--349.

\bibitem{COS10} C. Conti and L. Romani, A new family of interpolatory non-stationary subdivision schemes
for curve design in geometric modeling,in: \emph{ Numerical
Analysis and Applied Mathematics, International Conference} Vol-I
(2010).


\bibitem{COS11} C. Conti and L. Romani, Algebraic conditions on non-stationary subdivision symbols
 for exponential polynomial reproduction, \emph{J. Comput. Appl. Math.} , \textbf{236}
(2011), 543--556.



\bibitem{SDPS09} S. Daniel and P. Shunmugaraj, An approximating $C^2$
non-stationary subdivision scheme, \emph{Comput. Aided. Geom. D}
, \textbf{26} (2009), 810--821.
%

\bibitem{ND87} N. Dyn, J.A. Gregory and D. Levin, A 4-points interpolatory subdivision scheme
for curve design, \emph{Comput. Aided. Geom. D.} , \textbf{4(4)} (1987),
257--268.

\bibitem{ND95} N. Dyn and D. Levin, Analysis of asymptotically equivalent
binary subdivision schemes, \emph{J. Math. Anal. Appl.} , \textbf{193} (1995), 594--621.


\bibitem{SSK10} S.S. Siddiqi and N. Ahmad, An approximating $C^4$ stationary subdivision scheme,
 \emph{Eur. J. Sci. Res.}, \textbf{15(1)} (2006), 97--102.

\bibitem{ND02} N. Dyn and D. Levin, Subdivision schemes in geometric
modeling, \emph{ Acta Numerica} , \textbf{11} (2002), 73--144.




\bibitem{KC95} P.E. Koch, T. Lyche, M.
Neamtu and L. Schumker, Control curves and knot insertion for
trignometric splines, \emph{Adv. Comput. Math.} , \textbf{3}
(1995), 405--424.


\bibitem{JP12} J. Pana, S. Lin and X. Luo, A combined approximating and interpolating subdivision scheme with
$C^2$ continuity, \emph{Appl. Math. Lett.}, \textbf{25(12)}
(2012), 2140--2146.

\bibitem{SSK09} S.S. Siddiqi and N. Ahmad, A $C^6$ approximating subdivision scheme, \emph{Appl. Math. Lett.}, \textbf{21(7)} (2008), 722--728.


\bibitem{SS11} S.S. Siddiqi and N. Ahmad, A new three-point approximating $C^2$ subdivision
scheme, \emph{Appl. Math. Lett.}, \textbf{20(6)} (2007), 707--711.

\bibitem{SY12} S.S. Siddiqi and M. Younis, Construction of m-point binary approximating subdivision schemes, \emph{Appl. Math. Lett.} (2012), doi:10.1016/j.aml.2012.09.016.

\bibitem{SSN08} S.S. Siddiqi and N. Ahmad, A new five-point approximating subdivision scheme, \emph{Int. J. Comput. Math.}, \textbf{85(1)} (2008), 65--72.

\end{thebibliography}

\begin{figure}[tbp]
\begin{center}
 \includegraphics[width=2.7in, height=1.9in]{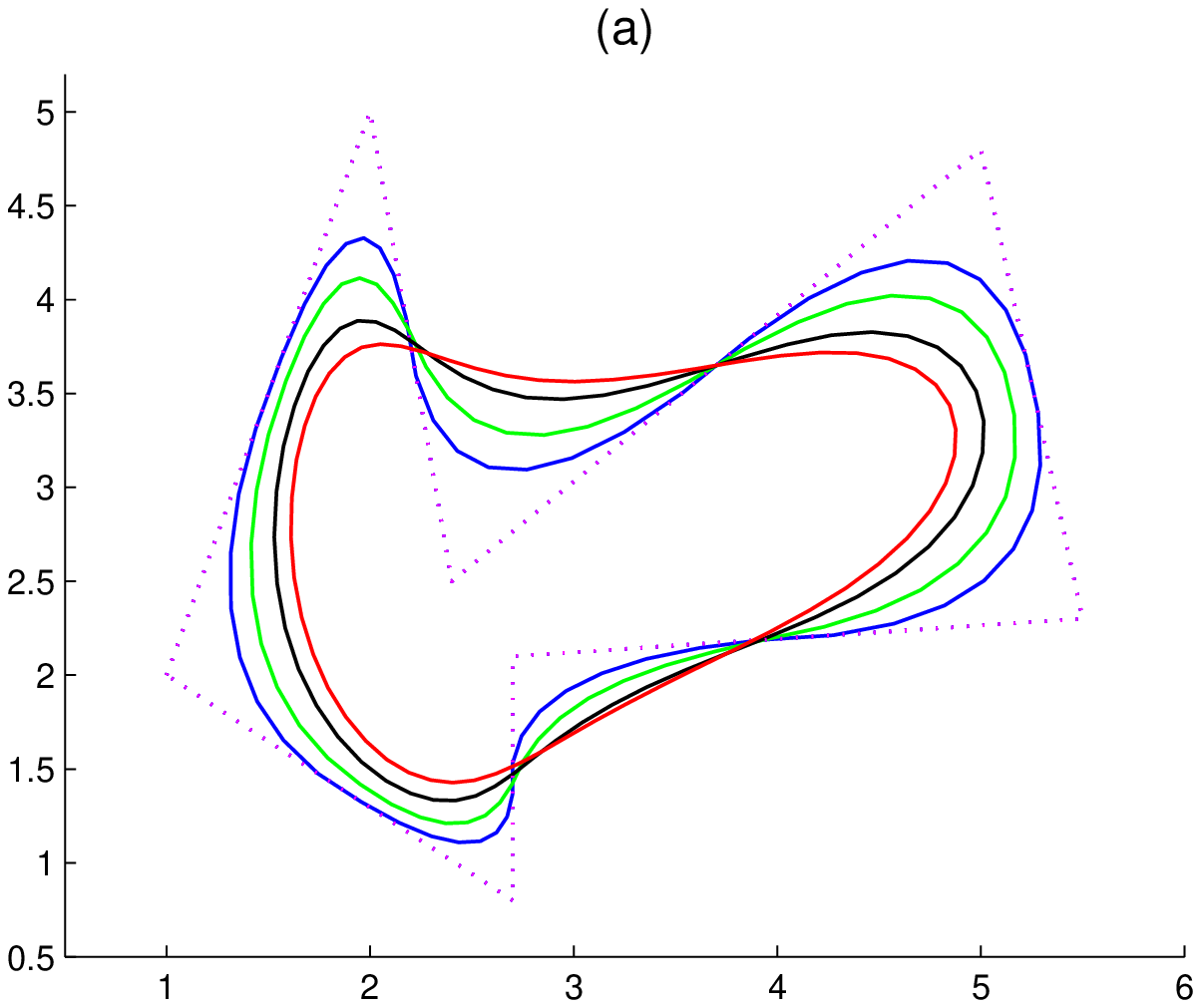}
  \includegraphics[width=2.7in, height=1.9in]{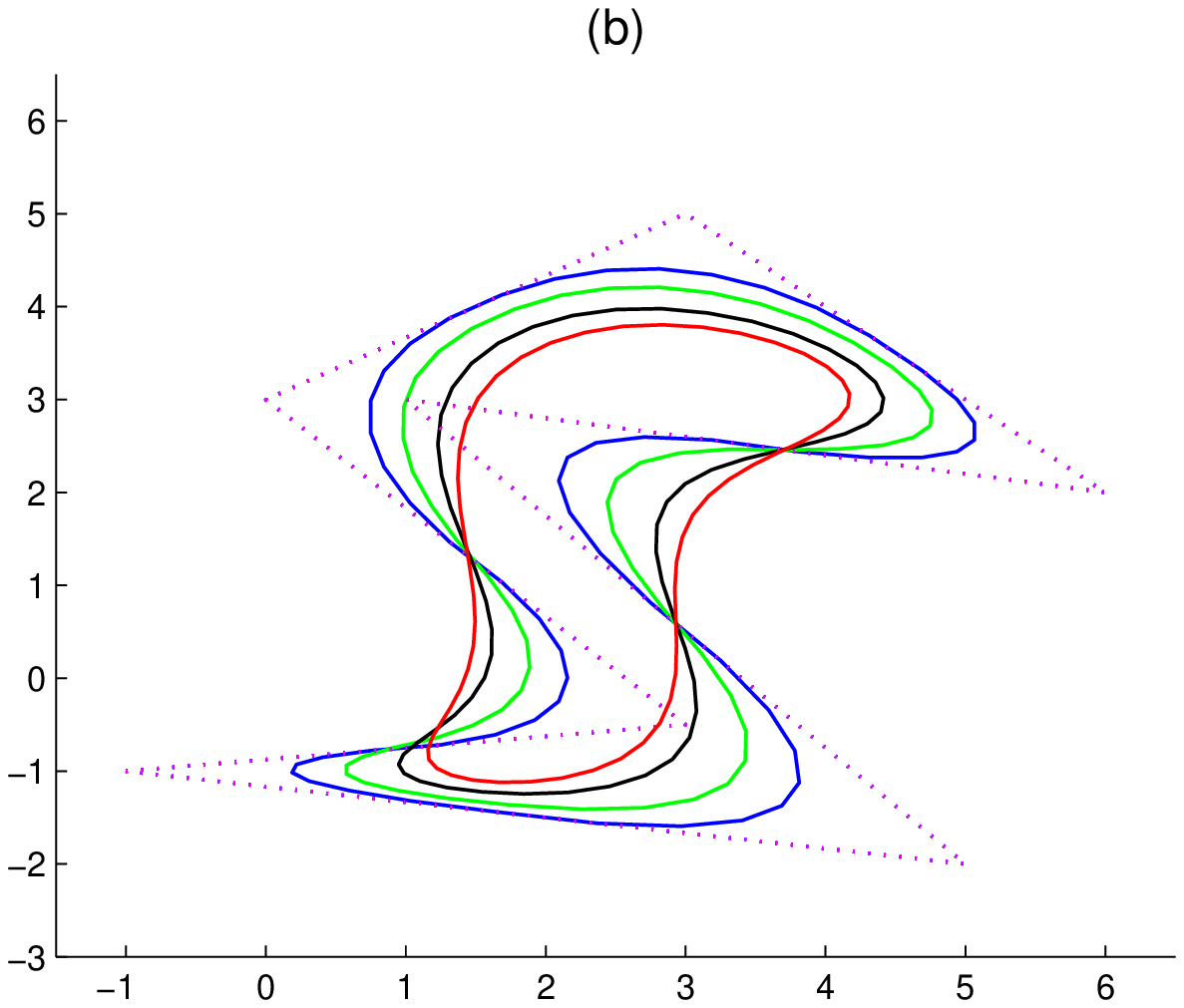}
   \includegraphics[width=2.7in, height=1.9in]{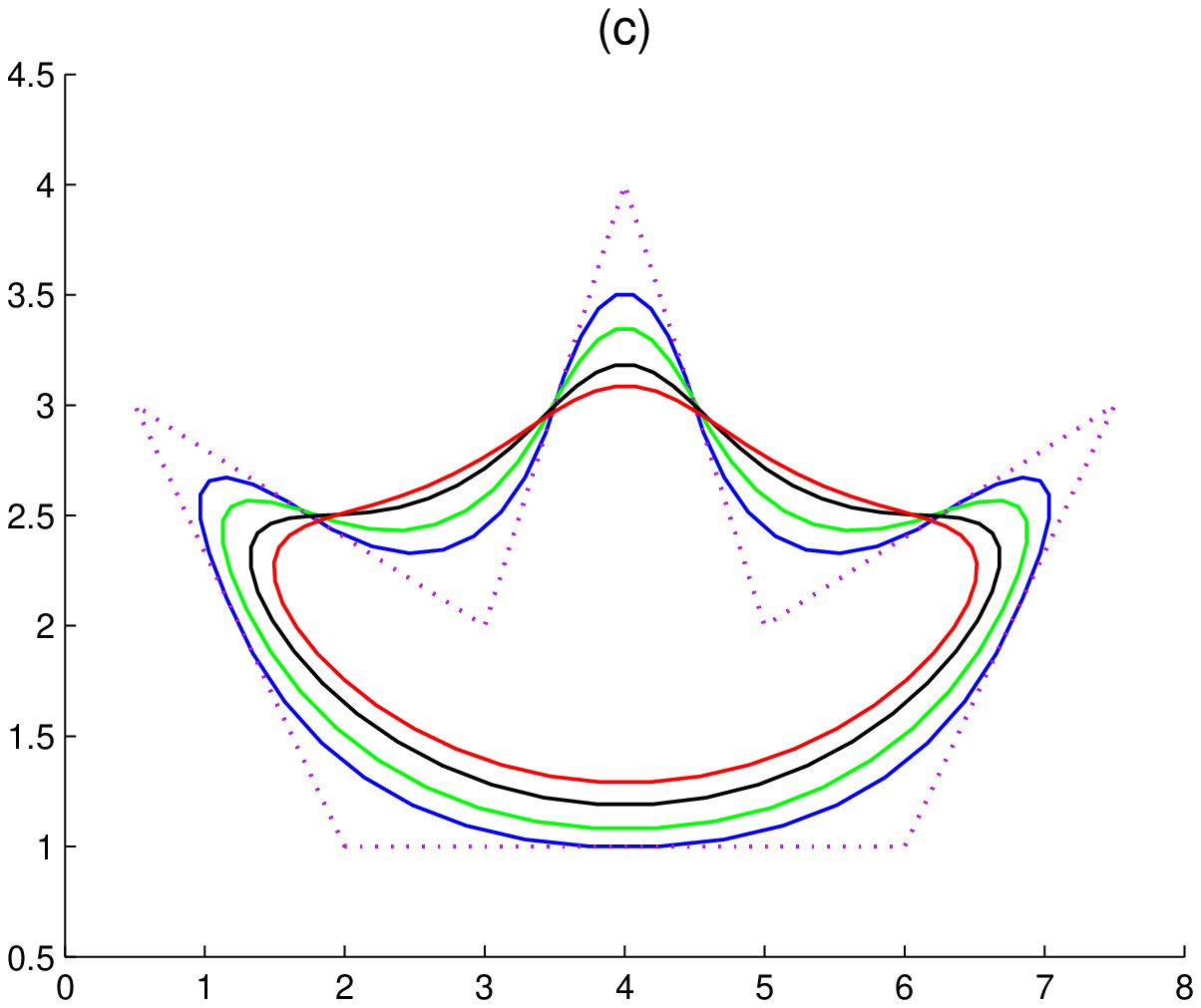}
   \includegraphics[width=2.7in, height=1.9in]{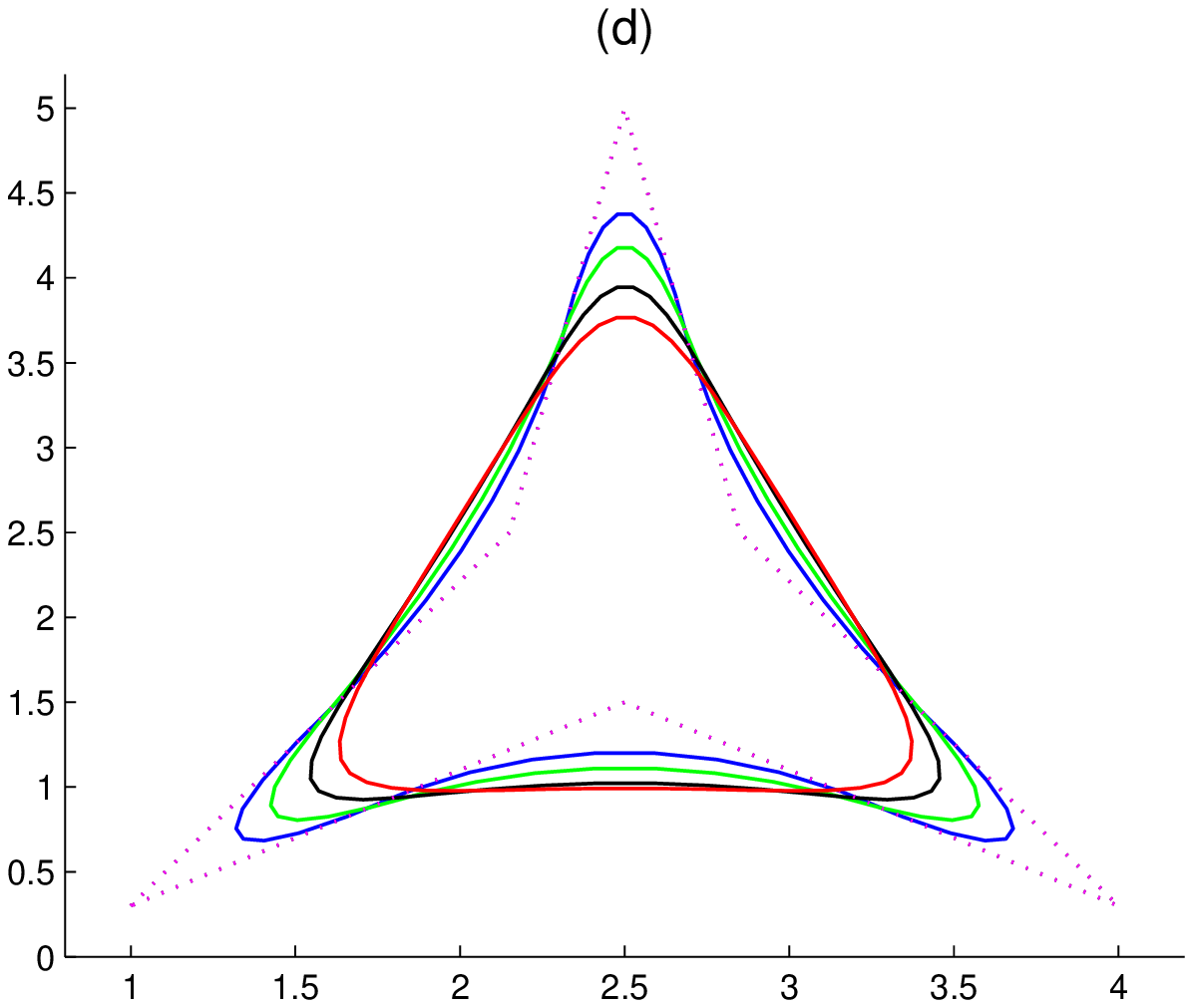}
 \caption{The continuous lines show the behavior of limit curves (after three
subdivision steps) of
the proposed scheme in blue, green, black and red colours,
respectively, after setting $m= 2, 3, 4$ and $5$.
The blue, green, black and red coloured limit curves coincide with
the schemes presented in \cite{GM74}, \cite{SDPS09}, \cite{SSK10}
and \cite{SSN08}, respectively.
 The dotted lines represent the control polygon}\label{fig:zzz}
\end{center}
\end{figure}

\begin{figure}[tbp]
\begin{center}
   \includegraphics[width=2.7in, height=1.9in]{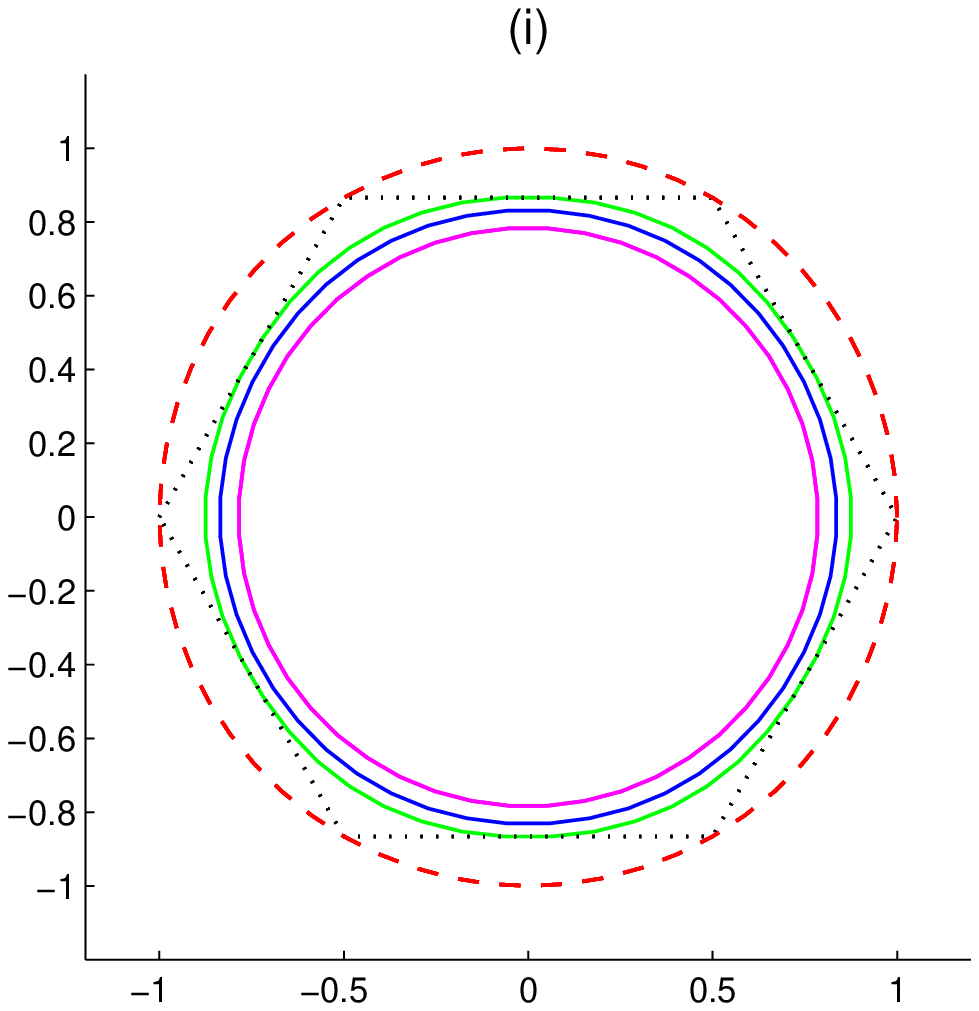}
   \includegraphics[width=2.7in, height=1.9in]{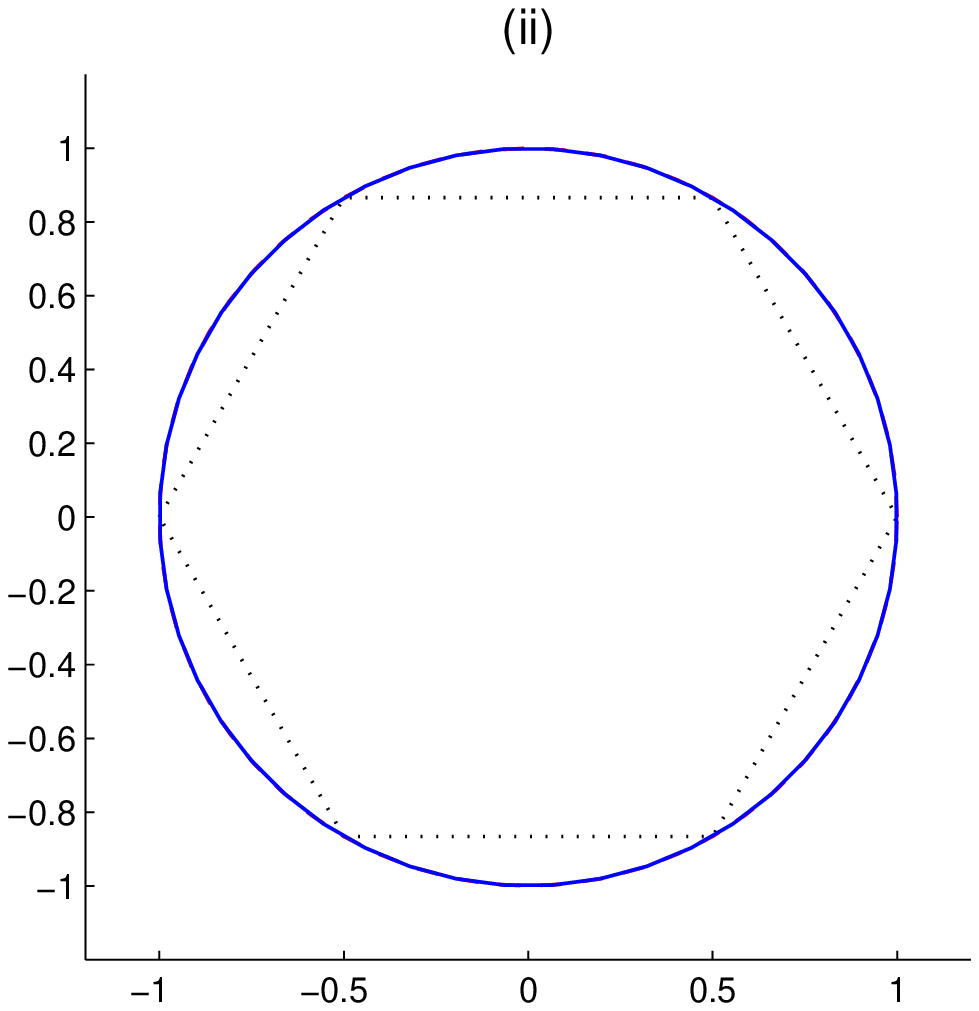}
\includegraphics[width=2.7in, height=1.9in]{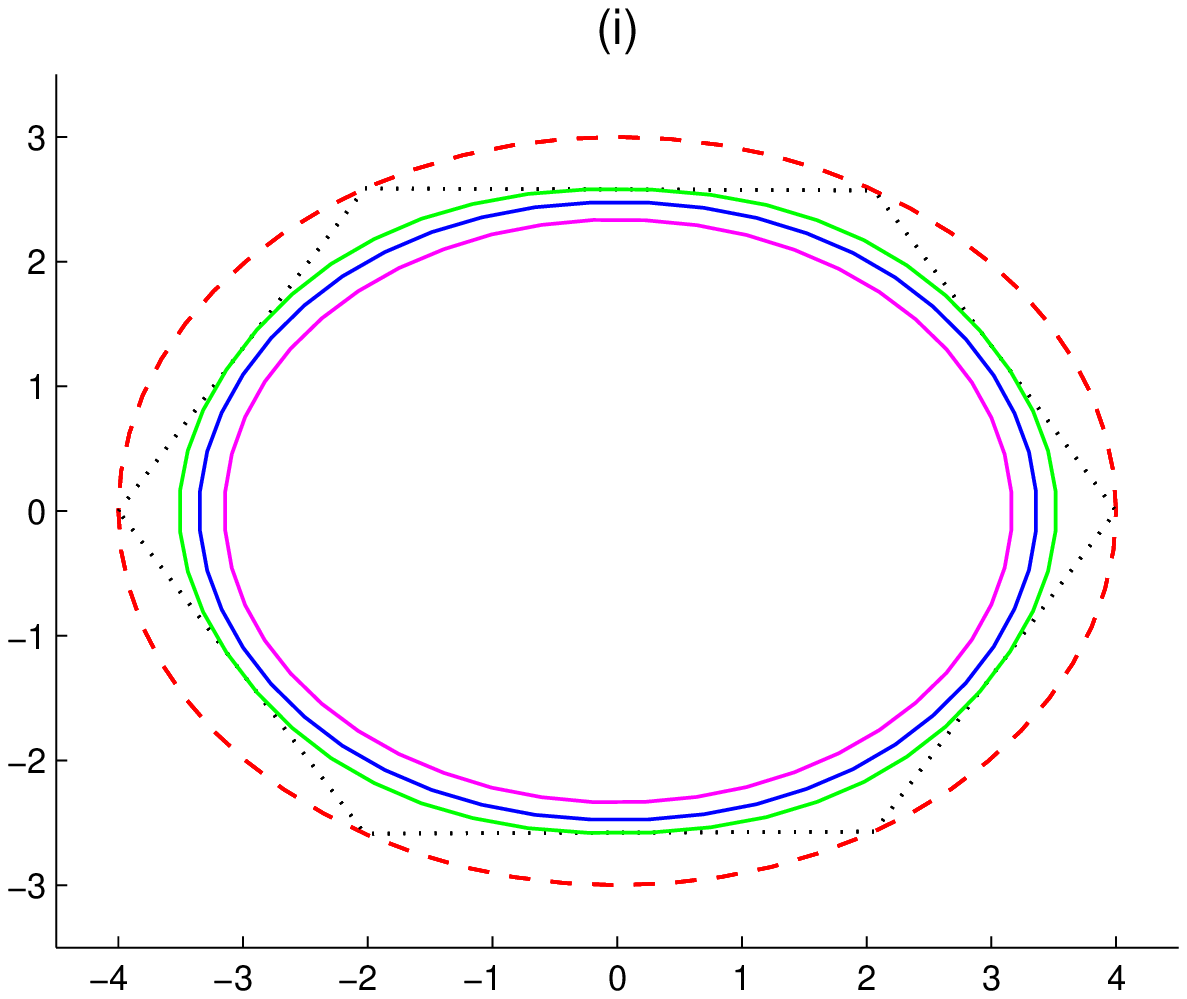}
\includegraphics[width=2.7in, height=1.9in]{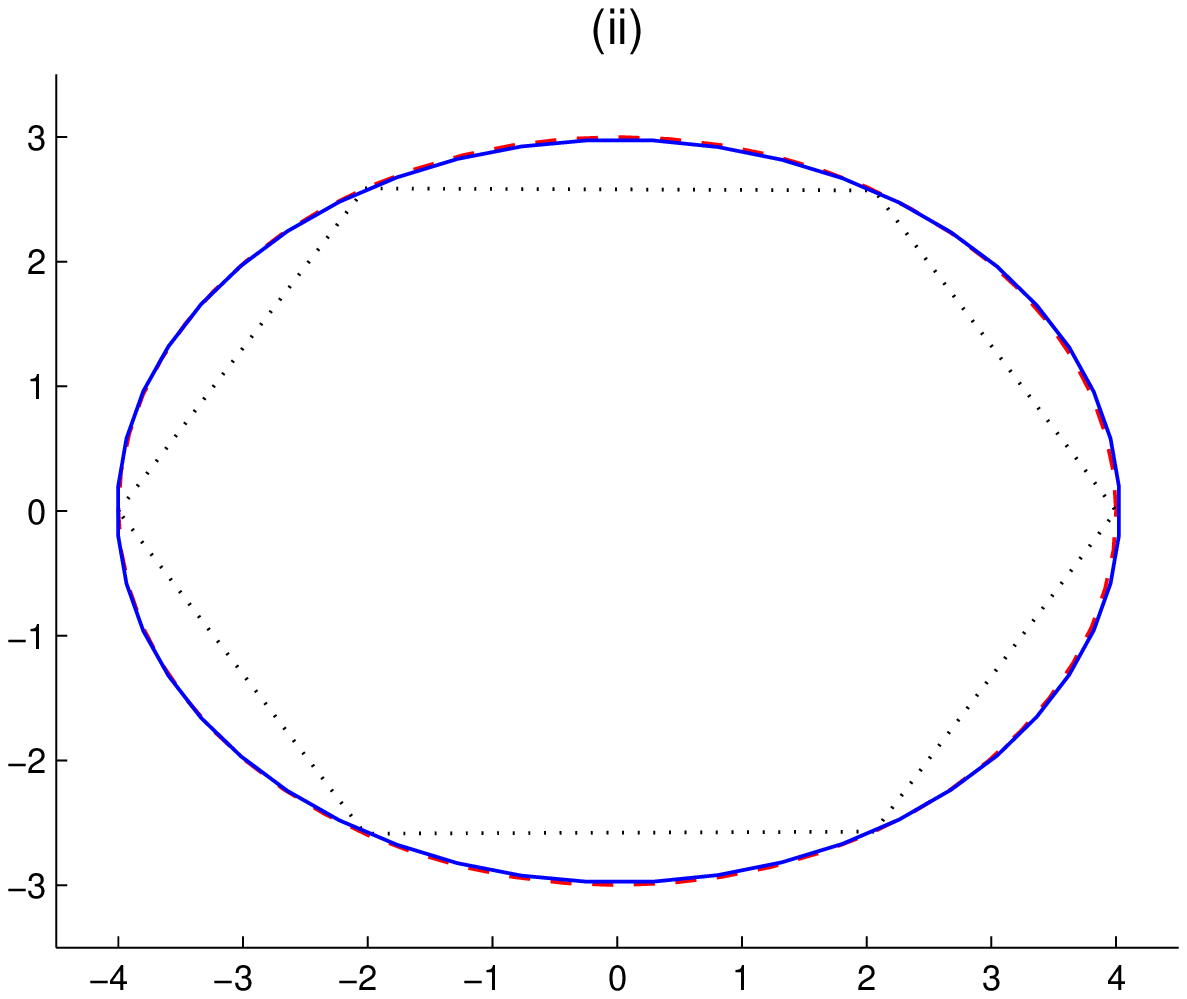}
\caption{The green, blue and magenta coloured continuous lines represent the limit curves of scheme \cite{GM74}, \cite{SS11} and \cite{SSK10}, respectively in (i). The blue coloured continuous line represents the limit curve of proposed scheme in (ii).The reproduction of the unit circle and ellipse (the limit curves in continuous lines are obtained after three subdivision steps). The red coloured broken lines represent unit circles and ellipses and black coloured dotted lines represent control polygon.}
\end{center}
\end{figure}

%
%

\end{document}